\documentclass[12pt,a4paper,oneside]{article}
\usepackage{amsmath,amsfonts,amssymb}
\usepackage{amssymb}
 \usepackage{amsthm}
\usepackage{amsmath}
 \usepackage{algorithmic}
 \usepackage{algorithm}
\usepackage[justification=centering]{caption}
\usepackage{tikz}
\usepackage{cite}
\usepackage{ulem}
\usepackage{url}
\usepackage{hyperref}
\usepackage{breakurl}
\usepackage{bigints}

\newlength{\defbaselineskip}
\newcommand{\setlinespacing}[1]%
           {\setlength{\baselineskip}{#1 \defbaselineskip}}

\newtheorem{theorem}{Theorem}[section]

\newtheorem{proposition}{Proposition}

\newtheorem{example}{Example}[section]

\title{{\bf\Large Inverse Fractional Knapsack Problem with Profits and Costs Modification}
{\footnotetext{
\hspace*{-0.7cm} N. T. Kien\\
 Department of Mathematics,
Teacher College, Cantho University, Cantho, Vietnam\\
e-mail: trungkien@ctu.edu.vn\vspace*{0.5cm}
\\
H. D. Quoc\\ Department of Mathematics, College of Natural Science,
Cantho University, Cantho,
Vietnam\\
e-mail: hdquoc@ctu.edu.vn }
}\\
\vspace{1cm}{\rm {\large{\bf Kien Trung Nguyen}
$\cdot$
 {\bf Huynh Duc Quoc}}} }
\date{}

\begin{document}

\maketitle

We address in this paper the problem of modifying both profits and costs of a fractional knapsack problem optimally such that a prespectified solution 
becomes an optimal solution with prespect to new parameters. This problem is called the inverse fractional knapsack problem. Concerning the $l_1$-norm, 
we first prove that the problem is $NP$-hard. The problem can be however solved in quadratic time if we only modify profit parameters. Additionally, 
we develop a quadratic-time algorithm that solves the inverse fractional knapsack problem under $l_{\infty}$-norm.

\section{Introduction}
The $\{0,1\}$ Knapsack Problem plays an important role in real-life decision making; for instance, finding the least wasteful way to cut raw materials, 
selection of investments and portfolios, cargo loading, etc. In order to solve this problem, we apply some non-polynomial but effective algorithms such as dynamic 
programming, greedy algorithm, branch and bound, and so on. The relaxation version of this problem is called the fractional knapsack problem, which can 
be solved by the greedy algorithm in $O(n\log n)$ time or by the algorithm of Balas and Zemel \cite{balas1980algorithm} in linear time. Here, we denote the input size 
of the problem by $n$.

The inverse (combinatorial) optimization problem consists of changing parameters of the problem at minimum total cost such that a prespecified solution 
becomes optimal with respect to new parameters. The first who investigated the inverse optimization problem were Burton and Toint \cite{Burton1992}. 
They developed an efficient algorithm that solves the inverse shortest path problem, which could be applied to predict the path of an earthquake. 
From here on, inverse optimization problem has increased interests from the community because of its potential applications. Ajuha et al. \cite{Ahuja1993} showed 
that the inverse linear programming optimization problem can be reduced to a problem of the same type, based on the so-called complementary slackness condition. 
In 2002, Ahuja and Orlin \cite{Ahuja1993} examined the inverse network flow problem with $l_1$- and $l_{\infty}$-norm. They presented combinatorial algorithms for solving 
this problem. Also, researchers focused on the inverse version of minimum spanning tree problem. Zhang et al. \cite{Zhang1997} was the first who investigated the problem with partition constraints
 in 1996 with practical applications.
Then Sokkalingam et al. \cite{Sokkalingam} solved the problem 
in $O(n^3)$ time. Ahuja and Orlin further improved the complexity of this problem to $O(n^2\log n)$. 
For terminology  concerning the inverse optimization problem and solution methods, readers refer to the survey of Heuberger \cite{Heuberger} .

Recently, the inverse $\{0,1\}$ knapsack problem has been investigated by Roland \cite{Roland}. He first 
showed that this problem under $l_{\infty}$-norm is co-$NP$-complete. Hence, there exists no approach to solve the problem in polynomial time, unless $P = NP$.
He also developed a pseudo-polynomial time algorithm based on a binary search to deal with the uniform-cost inverse $\{0,1\}$ knapsack problem. 
Besides, the author proposed a bilevel programming model for the problem under $l_1$-norm. Computation showed that this model is efficient enough.
In this paper we study the inverse fractional knapsack problem. According to the best of our knowledge, this problem has not been studied so far.

This paper is stated as follows. Section \ref{sec1} includes preliminary concepts and optimality criterion of the fractional knapsack problem. 
In Section \ref{sec2}, we formulate the inverse fractional knapsack problem under $l_1$-norm and show the $NP$-hardness in general case. If the cost coefficients
are fixed, a quadratic algorithm is developed. Section \ref{sec3} considers the inverse problem under $l_\infty$-norm. It is shown that the problem is solvable in $O(n^2)$ time.

\section{Problem definition}\label{sec1}
Let us first revisit the 0-1 knapsack problem and its fractional version. The 0-1 knapsack problem can be roughly stated as follows.
 Given a set of items, each with a cost and a profit, we wish to determine the number of each item to include in a collection so that the total cost is less than 
or equal to a given budget and the total profit is as large as possible. The relaxation of 0-1 Knapsack problem is the so-called fractional knapsack problem 
(FKP), which is formulated as follows.
\begin{equation}\label{EqK}
\begin{split}
\max\hspace{0.7cm}\displaystyle\sum_{i=1}^n p_{i}&x_{i} \\
\text{s.t.}\hspace{0.5cm} \displaystyle\sum_{i=1}^n c_{i}x_{i}&\leq b \\
0\leq x_{i}&\leq 1\quad\forall i=1,\dots, n.
\end{split}
\end{equation}

Here,the profits $p_{i}$ and the costs $c_i$ are positive intergers for all $i=1,\dots,n$. (FKP) can be solved by a simple greedy algorithm, where we takes 
the items with respect to the smaller ratios (profit over cost) until the budget constraint fulfills. Another solution approach with linear time complexity was proposed by 
Balas and Zemel \cite{Balas1980}, where we ruins a half of solution set until obtainning a stoping condition. For both of these two algorithms, an optimal 
solution has the following form.
\begin{proposition}
	Let $x^{0}=\left(x_{1}^{0},x_{2}^{0},\dots,x_{n}^{0}\right)$ be an optimal solution of \eqref{EqK}. Then all items in $x^0$ obtain value 0 or 1 except at most one item, say $x^{0}_{r}$, s.t. $0<x^{0}_{r}<1$.
\end{proposition}	
We now consider a solution $x^*$ where all items in $x^*$ obtain value 0 or 1. The conditions for $x^*$ to be an optimal solution of \eqref{EqK} can be straightforward derived from the greedy algorithm as below.
\begin{theorem}(Optimality Criterion)
	 A solution $x^*=\left(x^{*}_{1},x^{*}_{2},\dots,x^{*}_{n}\right)$ is an optimal solution of \eqref{EqK} iff the following conditions hold.
	 \begin{itemize}
	 	\item[\rm(i)]\rm $\sum_{i=1}^{n}c_{i}x_{i}^{*}=b$.
	 	\item[\rm(ii)]\rm $\displaystyle\min_{\{i:x^{*}_{i}\ne0\}}\left\{\dfrac{p_{i}}{c_i}\right\}\geq\displaystyle\max_{\{i:x^*_i=0\}}\left\{\dfrac{p_i}{c_i}\right \}$.

	 \end{itemize}
 
\end{theorem}


Next we formulate the inverse setting of (FKP). 
Given an instance of \eqref{EqK} and a prespecified solution $x^*$, where $x^*_i=0$ or $x^*_i=1$ for $i=1,\ldots ,n$. The profits and the costs can be either increased or reduced, i.e.,  $\tilde{p}_i=p_i+u_i-v_i$ and $\tilde{c}_i=c_i+\lambda_i-\mu_i $. Let $\left(u,v,\lambda,\mu\right)$ be the vector of modification and $C(u,v,\lambda,\mu)$ is the cost function. The inverse continuous knapsack problem is stated as follows:
\begin{itemize}
 \item The vector $x^*$ become an optimal solution of the modified knapsack instance.
 \item  The cost $C(u,v,\lambda,\mu)$ is minimized.
 \item $u_i$, $v_i$, $\lambda_i$, $\mu_i$ are feasible for all $ i=1,\dots,n $. It means
 $$u_i\in\left[0,\bar{u}_i\right]\cap\mathbb{Z}, v_i\in\left[0,\bar{v}_i\right]\cap\mathbb{Z}, \lambda_i\in\left[0,\bar{\lambda}_i\right]\cap\mathbb{Z}, \mu_i\in\left[0,\bar{\mu}_i\right]\cap\mathbb{Z}$$
\end{itemize}
Note that the modifications must be intergers to guarentee that the profits and costs are intergers, too. We can formulate the problem as a programming as below.

\begin{equation}\label{eq:}
\begin{split}
 \min \hspace{0.8cm} &C(u,v,\lambda,\mu) \\
\text{s.t.} \hspace{0.5cm} x^*&\in \textbf{argmax}\left\{\sum_{i=1}^{n}\tilde{p}_{i}x_i | \sum_{i=1}^{n}\tilde{c}_{i}x_{i}\leq b, x_i\in [0,1] \right\}\\
\tilde{p}_i &= p_i + u_i - v_i, \forall i = 1,\ldots,n \\
\tilde{c}_i &= c_i + \lambda_i - \mu_i, \forall i = 1,\ldots,n \\
u_i &\in [0,\bar{u}_i] \cap \mathbb Z, \forall i = 1,\ldots,n \\
v_i &\in [0,\bar{v}_i] \cap \mathbb Z, \forall i = 1,\ldots,n \\
\lambda_i &\in [0,\bar{\lambda}_i] \cap \mathbb Z, \forall i = 1,\ldots,n \\
\mu_i &\in [0,\bar{\mu}_i] \cap \mathbb Z, \forall i = 1,\ldots,n. 
\end{split}
\end{equation}
Here, objective function $C(u,v,\lambda,\mu)$ is a nondecreasing function. Recently,  one often considers the objectives w.r.t.
 $l_1$-, $l_2$-, $l_\infty$-norm, or Hamming distance for measuring the paying costs of the inverse optimization problems. In the following, we 
investigate properties and algorithms regarding (FIFKP) under $l_1$- and $l_\infty$-norm.

\section{The problem under $l_1$-norm}\label{sec2}
Let us consider (IFKP) under $l_1$-norm. Assume that we pay $w_i$ ($w'_i$) for modifying one unit of profit (cost),
 the corresponding objective function can be written as 
\begin{center}
$\displaystyle C(u,v,\lambda,\mu) := \sum_{i=1}^n (w_i(u_i + v_i) + w'_i(\lambda_i + \mu_i))$.
\end{center}

We first get the following result concerning (IFKP) under $l_1$-norm.
\begin{theorem}
(IFKP) under $l_1$-norm with both variable profits and costs is NP-hard.
\end{theorem}	
\textbf{Proof.}
Consider an instance of Partition problem (PP). Given a set of integers 
$S = \left\{a_1,a_2,...,a_n\right\}$ such that $\sum_{i=1}^{n}a_i=2B$, where $B$ is a positive integer. Does there exist a subset $S'$ of $S$ such that 
$\sum_{a_i\in S'}a_i = B$? This problem is NP-complete; see Garey and Johnson \cite{Garey}.

The decision version of (IFKP) is stated as follows. Given an instance of the inverse fractional knapsack problem. 
Does there exist a modification of profits and 
costs such that a prespectified solution become optimal and the objective value is at most $C$?

Given an instance of (PP). We construct an instance of (IFKP) in polynomial time.
\begin{itemize}
	\item The profits are $p_i := 4a_i$ for $i=1,\ldots, n$ and $p_{n+1} = 4$.
	\item The costs are $c_i := 2a_i$ for $i=1,\ldots ,n$ and $c_{n+1}= 1$.
	\item Let $\bar{\mu}_i := a_i$ and $\bar{u}_i := 4a_i$ for $i=1,\dots,n$; $\bar{\lambda}_i = \bar{v}_i := 0$ for $i=1,\dots, n+1$; 
$\bar{u}_{n+1} = \bar{v}_{n+1} := 0$.
	\item Set $b := 3B$ and $x^*=(1,\dots,1,0)$ with $n$ 1's and choose $C := 7B$.
\item The corresponding weight to modify one unit of profit $p_i$ is $w_i = 1$ and cost $c_i$ is $w'_i = 3$ for $i=1,\ldots,n$.
\end{itemize}
Observe that, in the current state of the problem we obtain
$\frac{p_i}{c_i} < \frac{p_{n+1}}{c_{n+1}}$ for $i = 1,\ldots ,n$. Furthermore, $\sum_{i=1}^{n+1}c_ix^*_i = 4B$.
Hence, vector $x^*$ is not feasible. To make it an optimal solution, we increase the profits $p_i$ for $i=1,\ldots, n$ and reduce the costs 
$c_i$ for $i=1,\ldots ,n$. In what follows we prove that the answer to (PP) is 'yes' iff the answer to (IFKP) is 'yes'.

Assume that the answer to (PP) is 'yes'. Then there exists a subset $S' \subset S$ such that $\sum_{a_i\in S'} = B$. 
We set $u_i = 0$ and $\mu_i := a_i$ for $a_i \in S'$. Otherwise, let $u_i := 4a_i$ and $\mu_i := 0$ for $a_i \not \in S'$. 
It is trivial to check that $x^*$ is an optimal solution of the modified fractional 
knapsack problem and the objective is $B$.

Conversely, assume that the answer to (IFKP) is 'yes'. We prove that the answer to (PP) is 'yes'. We first prove that the modification of 
profit $p_i$ can be shifted to the modification of $c_i$ without increasing the objective value. Indeed, assume that we modify profit $p_i$ by $x$ units and 
$c_i$ by $y$ units ($y \leq a_i$). By the optimality criterion, we get $\frac{4a_i+x}{2a_i-y} = 4$ or $x+4y = 4a_i$. Hence, $x$ is a multiplier of $4$.
 The modification yield an objective value $x + 3y = 4a_i - y \geq 3a_i$.
Furthermore, if we shift $\frac{x}{4}$ units from $p_i$ to $c_i$, we get $\frac{4a_i}{2a_i - \frac{x}{4} - y} = \frac{4a_i}{a_i} = 4$ and the corresponding 
objective value $3(\frac{x}{4} + y) = \frac{3}{4}(x+4y) = 3a_i$. Hence, it is trivial that the second option of modification reduce the objective value.
 We can assume that there exists at most one modification $\lambda_{i_0}$ such that 
$\lambda_{i_0} = a_{i_0} - k$ with $0 < k < a_{i_0}$ and $\lambda_i = 0$ or $\mu_i = a_i$ for $i \neq i_0$. 
Let us set $I := \{i \in \{1,\ldots,n\}: \mu_i =  a_i\}$. As $x^*$ become an optimal  solution of (FKP), we get  
\begin{center}
$ \displaystyle \sum_{i \in I} a_i + 2\sum_{j\not \in I}a_j - k = 3B $ or $ \displaystyle \sum_{i \not \in I}a_i = B + k$ \hfill $(1)$.
\end{center}
As the objective function is at most $B$, we get 
\begin{center}
$ 3(\displaystyle \sum_{i\in I} a_i + k) + 4\sum_{j\not\in I}a_j-2k \leq 7B$ or $ \displaystyle \sum_{j\not \in I} a_j - \frac{k}{2}  \leq B$ \hfill $(2)$.
\end{center}
From $(1)$ and $(2)$ we get $\frac{k}{2}\leq 0$. In other words, $k = 0$ and $\sum_{i\in I}a_i = B$. 
Set $S' := \{a_i: i\in I\}$, then the sum of items in $S'$ is $B$.

$\hfill\square$

As (IFKP) is $NP$-hard, there does not exist a polynomial time algorithm to solve it, unless $P = NP$. Therefore, approximation, heuristic approach, or 
special polynomially solvable cases of the problem are interesting topics. 

Now let $x^*$ be a feasible solution of (FKP), i.e., $\sum_{i=1}^n c_ix^*_i = b$.
We focus on the problem of modifying only profit parameters while the costs are fixed, i.e., we set $\lambda_i = \mu_i = 0$ for $i=1,\ldots,n$. We call  
this problem the fixed cost inverse fractional knapsack problem (FIFKP). It is trivial to get the following result.

\begin{proposition}\label{Prop1}
There exists an optimal modifications of (FIFKP)   s.t. the profits $\tilde{p}_i $ are increased if $i\in\left\{ i:x^*_i = 1\right\}$ and reduced if
 $i\in\left\{ i:x^*_{i}=0 \right\}$
\end{proposition}
\textbf{Proof.}
 It is straight forward as we have to increase the ratios $\frac{\tilde{p}_i}{c_i}$ for $i\in\left\{ i:x^*_i\ne 0\right\}$ and reduce the ones 
for  $i\in\left\{ i:x^*_{i}=0 \right\}$. \hfill $\square$

We denote by $I^0=\left\{i:x^*_i=0\right\}$ and $I^1=\left\{i:x^*_i = 1\right\}$.
We first consider the feasibility condition of (FIFKP). Let 
$$L=\displaystyle\max_{i\in I^0}\left\{ \dfrac{p_i-\bar{z}_i}{c_i} \right\},\quad U=\displaystyle\min_{i\in I^1}\left\{\dfrac{p_i+\bar{z}_i}{c_i}\right\}.$$
Then we obtain the following result.
\begin{proposition}
(FICKP) is feasible iff $L\leq U$.
\end{proposition}
{\bf{Proof.}}
$L$ is the maximum reduction of ratios $\frac{\tilde{p}_i}{c_i}$ for $i \in I^0$, while $U$ is the maximum augmentation of ratios 
$\frac{\tilde{p}_i}{c_i}$ for $i \in I^1$. The result follows. \hfill $\square$

From now on, we always assume that (FIFKP) is feasible.
By Proposition \ref{Prop1}, we set $v_{i}=0$ for $i\in I^1$ and $u_{i}=0$ for $i\in I^0$. 
We further denote by 
$z_i:= \begin{cases}
                u_i, & \mbox{if \,}  i\in I^1,\\
                 v_i,& \mbox{if \,}  i\in I^0,\\
 \end{cases}$
and
$\bar{z}_i:= \begin{cases}
                \bar{u}_i, & \mbox{if \,}  i\in I^1,\\
                 \bar{v}_i,& \mbox{if \,}  i\in I^0.\\
 \end{cases}$.
We say that a profit $p_i$ is modified by $z_i$ if it is reduced (increased) by $z_i$ for $i\in I^0$ ($I^1$).
Assume that modifying an amount of $p_i$ yields a corresponding cost $w_i$. 
The objective function can be written as
$$C(z)=\displaystyle\sum_{i=1}^{n}w_{i}z_i.$$

{\bf{Presolution:}} 
We presolve the problem by increasing the profits w.r.t. the items in $\{i\in I^1:\frac{p_i}{c_i} < L \}$ to $L$ and 
reducing the items in $\{i\in I^0:\frac{p_i}{c_i} > U \}$ to $U$. In other words,
For $ i\in I^1$ and $\dfrac{p_i}{c_i} < L$, we find the minimum value $z^0_i$  such that 
$\dfrac{p_i + z^0_i}{c_i} \geq L$. As $z_i$ is an integer, we can set $z^0_i := \lceil c_iL - p_i \rceil$. By the same argument, we set 
$z^0_i := \lceil p_i - c_iU \rceil$ for $i \in I^0$ and $\dfrac{p_i}{c_i} > U$.
The corresponding cost is $C_0 := \sum_{i\in \mathcal P}w_iz^0_i$, where $\mathcal P := \{i \in I^1:  \dfrac{p_i}{c_i} < L \} \cup 
\{i \in I^0:  \dfrac{p_i}{c_i} >U \}$.

Next we solve (FIFKP). Let $\alpha :=\displaystyle\min_{i\in I^1}\left\{\frac{p_i}{c_i}\right\}$ and 
$ \beta:=\displaystyle\max_{i\in I^0}\left\{\frac{p_i}{c_i}\right\}$. Also,  we denote by
$\tilde{I}^{\star} := I^{\star} \cap \left \{ i: \frac{p_i}{c_i}\in [\alpha,\beta] \right \}$ for $\star = 0,1$.
Observe that, we only modify the profits with indices in $\tilde{I}^0 \cup \tilde{I}^1$.
For a parameter $t \in [\alpha,\beta]$, we reduce (increase) the ratio $\frac{\tilde{p}_i}{c_i}$ such that it is less than (greater than) $t$ for 
$j \in \tilde{I}^0$ ($i \in \tilde{I}^0$). Denote by 
$I^1(t) := \{i \in \tilde{I}^{1}: \frac{p_i}{c_i } < t \}$ and $I^0(t) := \{i \in \tilde{I}^{0}: \frac{p_i}{c_i } > t \}$.
We find $z_i$ such that $\frac{p_i  - z_i}{c_i} \leq t$ or $z_i \geq p_i - c_it$ for $i \in I^0(t)$. As 
$z_i$ is an integer for $i \in I^0(t)$, we get $z_i := \lceil p_i - c_it \rceil$. Analogously, we can set $z_i := \lceil c_it - p_i \rceil$ for $i \in I^1(t)$.
Therefore, the objective function w.r.t. parameter $t$ can be written as follows.
$$C(t)=\displaystyle\sum_{i\in I^1(t)}w_i\lceil c_i t-p_i \rceil+\displaystyle\sum_{i\in I^0(t)}w_i \lceil p_i-c_{i}t \rceil .$$

Let $\mathcal B := \{\frac{p_i}{c_i}: i \in \mathcal P\}$. Assume that $\mathcal B := \{t_1,t_2,\ldots ,t_n\}$ with $t_1 < t_2 < \ldots < t_n$.
For $t \in (t_i,t_{i+1})$ with $t_i$ and $t_{i+1}$ being two consecutive members of $\mathcal B$, the set $I^0(t)$ and $I^1(t)$ do not change. 
In other words, $I^{\star}(t) = I^{\star}(t')$ for $\star = 0,1$ and $t,t' \in (t_i,t_{i+1})$. As $\lceil . \rceil$ is a quasi-concave function, we get the following result.
\begin{proposition}
$C(t)$ is a quasi-concave function for $t \in (t_i,t_{i+1})$, $t_i, t_{i+1} \in \mathcal B$.
\end{proposition} 

As $C(t)$ is quasi-concave for $t \in (t_i,t_{i+1})$, it is also quasi-concave in $[t_i,t_{i+1}]$. Therefore, the minimum value of $C(t)$ on $[t_i,t_{i+1}]$ 
is obtained at $t_i$ or $t_{i+1}$ as $C(t) \geq \min\{C(t_i),C(t_{i+1})\}$ for $t \in [t_i,t_{i+1}]$.

The following example states that $C(t)$ is however neither quasi-convex nor quasi-concave.

\begin{example}
	Given $x^*=(1,1,0,1,0)$ be a feasible solution. The corresponding profits and costs are given in the following table.

\begin{table}[!ht]
\centering
\renewcommand{\arraystretch}{1.25}
\begin{tabular}{|c|c|c|c|c|c|}
	\hline
	i&1&2&3&4&5\\
	\hline
	$p_i$&8&7&9&10&11\\
	\hline
	$c_i$&5&10&10&10&10\\
	\hline
	$\bar{z}_i$&3&4&3&1&4\\
	\hline
	$w_i$&3&$1/2$&$1/2$&1&1\\
	\hline
	\end{tabular}
\caption{An instance of (FICKP)}\label{tab1}
\end{table}
First of all, the set $\mathcal B$ consists of
$t_1=\dfrac{3}{5}; t_2=\dfrac{7}{10}; t_3=\dfrac{9}{10}; t_4=\dfrac{10}{10}; t_5=\dfrac{11}{10}$.
We compute the objective value at each break-points as follows.
$C(t_1)=\dfrac{3}{2}+5=6,5; C(t_2)=3+1+4=8; C(t_3)=6+1+2=9; C(t_4)=6+\dfrac{3}{2}+1=8,5; C(t_5)=9+2+1=12$.
Hence, $C(t)$ is neither quasi-convex nor quasi-concave for $t\in [t_1=\dfrac{3}{5};t_1=\dfrac{11}{10}]$.
\end{example}

Now we know that the objective function is neither quasi-convex nor quasi-concave. To find the optimal solution of $C(t)$, we first compute the values 
of $C(t)$ at all break-points in $\mathcal B$. Then we take the best one. The value of $C(t)$ at each break-point can be computed in linear time. Furthermore,
there are at most linearly many break-points. Hence, the optimal solution of $C(t)$ can be found in quadratic time. We get the main result of this section.
\begin{theorem}
The inverse fractional knapsack problem with variable profits can be solved in quadratic time.
\end{theorem}


\section{Problem under $l_{\infty}$-norm}\label{sec3}
Now we investigate the uniform-cost inverse fractional knapsack problem under $l_{\infty}$-norm. The corresponding objective function can be rewritten as follows.
\begin{center}
$\displaystyle \max_{i=1}^n\left\{u_i,v_i,\lambda_i,\mu_i\right\}$.
\end{center}

Let us recall that $I^1$ and $I^0$ are the set of items in $x^*$ with value 1's and 0's, respectively. Moreover, a property of modifying  profits and costs 
is given as follows.

\begin{proposition}\label{Prop2}
There exists an optimal solution such that we increase (reduce) the profits of items in $I^1$ ($I^0$) and increase the costs of items in $I^0$.
\end{proposition}
\textbf{Proof.}
Similar to Proposition \ref{Prop1}.\hfill$\square$

By Proposition \ref{Prop2}, we set $v_i := 0$ for $i\in I^1$, $u_i = \mu_i :=0$ for $i\in I^0$. We study the two following situations.

\textbf{Case 1:} If $\sum_{i=1}^n c_ix^*_i = b$, then  $x^*$ is a feasible solution. Therefore, we do not modify the cost coefficients $c_i$ for $i\in I^1$. 
Otherwise, if we do modity the costs of items in $I^1$, the optimality criterion does not hold according to the infeasibility of $x^*$. 

Let us set $\lambda_i =\mu_i :=0,\forall i\in I^1$.  By the optimality criterion, the following condition must hold
	$$\min_{i\in I^1}\left\{\dfrac{\tilde{p}_i}{\tilde{c}_i}\right\}\geq\max_{j\in I^0}\left\{\dfrac{\tilde{p}_j}{\tilde{c}_j}\right\}.$$
	Hence, for each $i, j$ such that $i\in I^1$, $j\in I^0$ and $\dfrac{p_i}{c_i}<\dfrac{p_j}{c_j}$, we calculate the minimum object value such that 
$$\dfrac{\tilde{p}_i}{c_i}\geq\dfrac{\tilde{p}_j}{\tilde{c}_j}$$
	Replacing $\tilde{p}_i,\tilde{p}$, and $\tilde{c}_j$ by $p_i+u_i,p_j-v_j$, and $c_j+\lambda_j$, we get
	$$\dfrac{p_i+u_i}{c_i}\geq\dfrac{p_j-v_i}{c_j+\lambda_j}$$
After some elementary computations, we get the inequality
\begin{equation}\label{eqRatio}
u_i\lambda_j+p_i\lambda_j+u_{i}c_j+c_{i}v_i\geq c_{i}p_j-p_{i}c_j 
\end{equation}

We sort the corresponding upper bounds $\bar{u}_i, \bar{\lambda}_j, \bar{u}_i, \bar{v}_i$, then we compute the value on the left hand side of \eqref{eqRatio} 
w.r.t. the threholds. Then it is trivial to compute the smallest objective such that \eqref{eqRatio} holds.
Let $K_{ij}$ be the minimum value such that $\dfrac{\tilde{p}_i}{c_i}\geq\dfrac{\tilde{p}_j}{\tilde{c}_j}$. Then 
$K=\displaystyle\max\min_{\substack{i\in I^1\\j\in I^0}}K_{ij}$ is the optimal object value.

\textbf{Case 2:} If $\sum_{i=1}^n c_ix^*_i < b$, the vector $x^*$ is not feasible. 
Therefore, we first modify the cost optimally so that $x^*$ become a feasible solution as follows.

Let $S := \{\bar{\mu_i}: i \in I^1\}$. We aim to find the smallest value, say $\bar{\mu}_{i^*}$, in $S$ such that 
\begin{center}
$\frac{b-\sum_{j\in I^1}c_j}{|I^1|} \leq \bar{\mu}_{i^*}$ 
\end{center}
for $i\in I^1$.
It can be done by applying a binary search algorithm. Indeed, let $m$ be the median of $S$ and $\bar{\mu}_{i_0}$ be the largest element in $S$ which is less than or equal $m$. 
If $\frac{b-\sum_{j\in I^1}c_j}{|I^1|} < \bar{\mu}_{i_0}$, we know that  $\mu_{i^*} \leq \mu_{i_0}$ and one thus has to find $\mu_{i^*}$ in 
$S := S\backslash \{\mu_i > \mu_{i_0}\}$. 
Otherwise, we know that $\mu_{i^*} > \mu_{i_0}$ and consider  $S := S\backslash \{\mu_i \leq \mu_{i_0}\}$. 
This algorithm find $\bar{\mu}_{i^*}$ in linear time.

After finding $\bar{\mu}_{i^*}$, we set $\mu_i:=\bar{\mu}_i$ for $i\in I^1$ and $\bar{\mu}_i < \bar{\mu}_{i^*}$. Then set 
$\mu_i := \max\{\mu_i: \mu_i \in S \text{ and } \bar{\mu}_i \geq \bar{\mu}_{i^*}\}$ and $J^1 := \{i \in I^1: \bar{\mu}_i \geq \bar{\mu}_{i^*}\}$
and $J^2 := \{i \in I^1: \bar{\mu}_i < \bar{\mu}_{i^*}\}$. 
We further find  $\mu^{\min} = \left \lceil  \frac{b - \sum_{i\in I^1}c_i - \sum_{j\in J^2}\bar{\mu}_j}{|J^1|} \right \rceil$. Then it is easy to verify that 
\begin{center}
$\displaystyle \sum_{j\in J^2}(c_j+\bar{\mu}_j) +  \sum_{j\in J^1}(c_j+\mu^{\min}) \geq b $.
\end{center}

Next we study which variables in $J^1$ should takes value $\mu^{\min}$ and $\mu^{\min}-1$. We first consider how many variables in $J^1$ 
take value $\mu^{\min}$. This number equals 
\begin{center}
$N := b - \displaystyle \sum_{j\in J^2}(c_j+\bar{\mu}_j) - \sum_{j\in J^1}(c_j+\mu^{\min} - 1)$.
\end{center}
Hence, there are $N$ variables in $J^1$ obtainning value $\lambda^{\min}$ and $|J^1| - N$ variables in $J^1$ obtainning value $\lambda^{\min} - 1$.

The objective function is $\lambda^{\min}$ to modify the costs so that $x^*$ becomes feasible. In order to reduce the ratios $\frac{\tilde{p}_i}{\tilde{c}_i}$
for $i \in I^0$ and augment the the ratios $\frac{\tilde{p}_i}{\tilde{c}_i}$ for $i \in I^1$, we set 
$\tilde{p}_i := \begin{cases}
                 p_i +  \min\{\bar{u}_i,\lambda^{\min}\}, & \mbox{if \,}  i\in I^1,\\
                  p_i -  \min\{\bar{v}_i,\lambda^{\min}\}, & \mbox{if \,}  i\in I^0,\\
 \end{cases}$.
We also update 
$\bar{\star}_i := \begin{cases}
                 0, & \mbox{if \,} \bar{\star}_i \leq \lambda^{\min},\\
                  \bar{\star}_i - \lambda^{\min} , & \mbox{if \,} \bar{\star}_i > \lambda^{\min}\\
 \end{cases}$ for $\star = u, v$ and $i\in I^0$ or $i\in I^1$, accordingly.

Let us now consider the set the current rations $\frac{\tilde{p}_i}{\tilde{c}_i}$. For each $i \in I^1$ we compute the largest cost, say $K_i$, such that 
 $\frac{p'_i}{c'_i} \geq \frac{p'_j}{c'_j}$ for all $j \in J^0$ by the similar approach in Case 1. We then consider the chance to reduce the cost as follows.
We first take $|J^1| - N$ items with respect to the $|J^1| - N$ largest costs in $\{K_i\}_{i\in I^1}$ and set $\tilde{p}_i = \tilde{p}_i + \lambda^{\min} - 1$. 
Then we reevaluate the costs $K'_i$ with respect to new items. 

\textbf{Case 3:} If $\displaystyle\sum_{i=1}^n c_ix^*_i > b$, we can solve the problem as in Case 2.

In summary, we first check the feasibility of $x^*$. If it is not feasible, we can justify the cost coefficients in linear time to make it feasible. Then, 
it costs quadratic time to compute $K_{ij}$ for $i\in I^1$ and $j \in I^0$. The final step costs also quadratic time in order to find $\tilde{K}_{ij}$. Hence, the total 
computation complexity is quadratic.

\begin{theorem}
The inverse fractional knapsack problem under $l_\infty$-norm can be solved in quadratic time.
\end{theorem}

\section{Conclusions}
We considered the inverse fractional knapsack problem with profit and cost modifications. It is shown that the problem under $l_1$-norm is, in general, $NP$-hard.
Especially, if we can only justify the profit parameters, this problem is solvable in $O(n^2)$ time. Moreover, we solve the problem under $l_\infty$-norm in quadratic time based on 
greedy type algorithm. 

It is promising to study the inverse fractional knapsack problem under various of objective functions. Furthermore, it is also worthwhile to consider the inverse 
mixed integer knapsack problem by combining the techniques in  this paper and in Roland \cite{Roland}.

\end{document}